\documentclass[11pt]{amsart}
\usepackage{amssymb}
\usepackage{amsmath}
\usepackage{amsfonts}
\usepackage{color}
\newtheorem{theorem}{Theorem}
\theoremstyle{plain}
 
\newtheorem {definition}[theorem]{Definition}

\newtheorem {prop}{Proposition}

\newtheorem {theo}{Theorem}
\newtheorem {corol}{Corollary}
\newtheorem {example}{Example}
\usepackage{soul}
\usepackage{ulem}
\usepackage{cancel}
\begin{document}
\renewcommand{\bibname}{Bibliography}
\title{On polynomially hypo-EP operators. }

\author{R. Semmami}
\address{R. Semmami,  Laboratory of Mathematical Analysis and Applications, Faculty of Sciences, Mohammed V University in Rabat, Rabat, Morocco.}
\email{semmami.rachi@gmail.com;rachid\_semmami@um5.ac.ma}
\author{H. Ezzahraoui}
\address{H. Ezzahraoui, Laboratory of Mathematical Analysis and Applications, Faculty of Sciences, Mohammed V University in Rabat, Rabat, Morocco.}
\email{hamid.ezzahraoui@fsr.um5.ac.ma;h.ezzahraoui@um5r.ac.ma}

\begin{abstract}
The purpose of this paper is to present a new class of operators known as polynomially hypo-EP operators,
extending the notation of hypo-EP, $n$-hypo-EP, and polynomially EP. The paper explores numerous properties and characterizations of polynomially hypo-EP operators and polynomially EP operators. In addition, we characterize polynomially hypo-EP operators applying the adequate operator matrix representations. By employing these findings, we derive new characterizations of EP, hypo-EP, $n$-EP and $n$-hypo-EP.
\end{abstract}

\maketitle

\section{Introduction}
Let $\mathcal{H}$ be a complex Hilbert space and let $\mathcal{L(H)}$ represent the space of
all bounded linear operators defined on $\mathcal{H}$ and $\mathcal{R(H)}\subset \mathcal{L(H)}$ be the set of  operators with closed range. For $T \in \mathcal{L(H)}$, $T^{*},$ $N(T)$ and $R(T)$ stand for its adjoint, kernel and range, respectively. For $T,S\in \mathcal{L(H)}$, we write $[T,S]:=TS-ST$ for the commutator of $S$ and $T.$

A closed subspace $E\subset\mathcal{H} $ is said to be $T$-\emph{invariant}  if $T(E) \subset E$. We define $T_{E}: E\rightarrow E$
as $T_{E}x=Tx,$ $x\in E.$ Also $E$ it is called  \emph{  reducing  } for $T$ if it is invariant for both $T$ and $T^*$, or equivalently, if $T(E)\subset E$ and $T(E^{\perp})\subset E^\perp$, where  $E^\perp$ is the orthogonal complement of $E$.

For   $T\in \mathcal{R(H)}$, the Moore-Penrose inverse $T^{\dagger}$ of $T,$ is defined as the unique solution of the
following four operator equations: $$TT^{\dagger}T=T,\;T^{\dagger}TT^{\dagger}=T^{\dagger},\;(T^{\dagger}T)^{*}=T^{\dagger}T\;\text{and }(TT^{\dagger})^{*}=TT^{\dagger}.$$
It follows immediately that
$$TT^{\dagger}=P_{R(T)},\; T^{\dagger}T=P_{R(T^{*})},\; N(T^{\dagger})=N(T^{*})\mbox{ and } N(T^{\dagger*})=N(T).$$
The Cauchy dual of an operator $T\in\mathcal{L(H)}$ with closed range is defined as $\omega(T)=T(T^{*}T)^{\dagger}.$ For  equivalent definitions and further properties of $\omega(T)$, see \cite{emz1}.

The concept of EP matrices was defined by Schwerdtfeger \cite{schwerdtfeger} and generalized to Hilbert space operators with closed range in \cite{camb}. An operator $T\in \mathcal{L(H)}$ is called EP if $R(T)$ is closed and $R(T)=R(T^{*})$. EP operators are significant because they commute with their Moore–Penrose inverse. Normal operators with a closed range and nonsingular operators are particular cases of EP operators. Several authors have thoroughly examined different established characterizations of EP matrices and EP operators. See for example \cite{Bourdon, camb,chen,driv, emz3,ferr, dija}.

An operator $T\in\mathcal{R(H)}$ is SD (star-dagger), if  $[T^*,T^\dag]=0$. This class been studied in detail by several authors. See  \cite{emz3,ferr,hartwig}.

The class of $n$-EP matrices was introduced in \cite{SLN} as a generalization of
EP matrices. The notion of $n$-EP operators was defined in \cite{wangchun} extending the
definition of $n$-EP matrices. An operator $T\in \mathcal{L(H)}$ is called $n$-EP and we write $T\in$ ($n$-EP) if $R(T)$ is closed and $T^{n}T^{\dagger}=T^{\dagger}T^{n}$ for some $n\in\mathbb{ N}$. For further properties of $n$-EP operators, we refer to \cite{emz1,menkad}.  Generalizing the notion of $n$-EP operators, polynomially EP operators were presented in \cite{dija}.

M. Itoh introduced in \cite{itoh} the class of hypo-EP operators   as follows. An operator $T\in \mathcal{R(H)}$ is an  hypo-EP operator if  $[T^{\dagger},T]\geq 0$. Let (HEP) be the class of  hypo-EP operators. Several authors have thoroughly examined different established characterizations of hypo-EP matrices and hypo-EP operators. See \cite{Bourdon, emz3, ding1, vinoth} for detailed study. Some research characterized hypo-EP operators through factorizations and  when the product of two hypo-EP matrices or hypo-EP operators is also hypo-EP, was studied. See for example \cite{Bourdon,joh, ding1}.

As an extension of hypo-EP operators, $n$-hypo-EP operators were introduced in \cite{emz3} as follows:
$T\in \mathcal{R(H)}$ is $n$-hypo-EP $(n\in \mathbb{N})$ if $R(T^{ n})\subseteq R(T^{*}).$  Let ($n$-HEP) be the class of $n$-hypo-EP operators.

  The operator $T$ is called normal if $[T^{*},T]=0.$ As an extension of normal operators, $n$-normal operators were introduced in \cite{jibril} as follows:
$T\in \mathcal{L(H)}$ is $n$-normal $(n\in \mathbb{N})$ if $[T^{ n}, T^{*}]=0$. Generalizing the
notion of $n$-normal operators, polynomially normal operators were presented
in \cite{djor}.

We use $Poly_{1}$ to denote the set of all complex polynomials in one variable.
Taking conjugate coefficients of $ p\in Poly_{1},$ we obtain $\overline{p}\in Poly_{1}$: $\overline{p}(z):= \overline{p(\overline{z})},$
 $z\in\mathbb{C}.$ Let $T\in \mathcal{L(H)}$ and $p\in Poly_{1}$ be nontrivial. If $[p(T),T^{*}]=0,$ then
we say that $T$ is $p$-normal. $T\in \mathcal{L(H)}$ is polynomially normal, if $T$ is $q$-normal for some $q\in Poly_{1}.$ Observe that, by \cite[Lemma 1]{djor}, $T\in \mathcal{L(H)}$ is $p$-normal if and only if $p(T)$ is normal. Recent results related to polynomially normal operators can be found in \cite{cvet}.

 Let $T\in \mathcal{L(H)}$ and $p\in Poly_{1}$ be nontrivial. Following \cite{dija}, we say that $T$ is $p$-EP, if $p(T)T^{\dagger}=T^{\dagger}p(T)$. An operator $T\in \mathcal{L(H)}$ is polynomially EP, if $T$ is $q$-EP for some $q\in Poly_{1}.$ It is proved in \cite[Lemma 2.4]{dija} that $T$ is $p$-EP if and only if $T^*$ is $\bar{p}$-EP where $(p(T))^*=\bar{p}(T^*)$.
%Recent results related to polynomially EP operators can be found in \cite{dija}.
It has been shown in  \cite{dija}  that every polynomially normal operator with a closed range is a polynomially EP operator and under additional assumptions. It has been verified  that  the converse is valid.

Motivated by recent research and useful properties of EP, hypo-EP, $n$-EP, $n$-hypo-EP, polynomially normal operators and polynomially EP operators, we define a new class of closed range operators called polynomially hypo-EP operators. The class of polynomially hypo-EP operators is a wider class than the classes of all hypo-EP, $n$-hypo-EP and polynomially EP operators. We prove many properties and characterizations
of polynomially hypo-EP operators and we give some additional properties of polynomially EP operators.
%Several of these properties are generalizations of well-known properties of EP, hypo-EP, $n$-EP and $n$-hypo-EP operators, and some of them are new characterizations of $n$-hypo-EP and $n$-EP operators.
We show that every polynomially EP operator is a polynomially hypo-EP operator and under additional assumptions, we verify that the converse is valid.  What is more, we characterize polynomially hypo-EP operators by the adequate operator matrix representations (Theorems \ref{rep} and  \ref{equi}) and we obtain new characterizations of $n$-hypo-EP and hypo-EP operators.

\section{Polynomially EP operators}
It’s easily seen that if $T$ is EP, then $T$ is $p$-EP for an arbitrary $p\in Poly_{1}.$ So, we proposed a wider class of operators which includes EP and $n$-EP operators. Obviously, when $T$ is $p$-EP for $p(t)=t$ (or $p(t)=t^{n}),t\in\mathbb{C},$ then
$T$ is EP ($n$-EP). From \cite[Exemple 2.2]{dija}, we observe that a polynomially EP operator is not EP in general. In the following proposition, we give a necessary condition of $p\in Poly_{1}$ such that an operator $m$-EP is $p$-EP.
\begin{prop}\label{croi}
  Let $T\in (m\text{-EP})$ where $m\geq1.$ Then $T\in (p\text{-EP})$ for all $p(t)=\sum_{j=0}^{n}a_{j}t^{j}\in Poly_{1}$ such that $m\leq k=\min\{j\geq1: a_{j}\neq0\}.$
\end{prop}

\begin{proof}
   Suppose that $T$ is $m$-EP and $p(t)=\sum_{j=0}^{n}a_{j}t^{j},$ where $a_{0},a_{1},\ldots,a_{n}\in\mathbb{C}.$

We have $$p(T)=a_{n}T^{n}+a_{n-1}T^{n-1}+\cdots+a_{k}T^{k}=T^{m}(a_{n}T^{n-m}+\cdots+a_{k}T^{k-m}).$$
Thus $$R(p(T))=R(T^{m}(a_{n}T^{n-m}+\cdots+a_{k}T^{k-m}))\subseteq R(T^{m})\subseteq R(T^{*}).$$
And  $$R(p(T)^{*})=R(T^{*m}(a_{n}T^{*n-m}+\cdots+a_{k}T^{*k-m}))\subseteq R(T^{*m})\subseteq R(T).$$
From \cite[Theorem 2.5]{dija}, $T$ is $p$-EP.
\end{proof}
In the special case that $p(t)=t^{n}$ in Proposition \ref{croi}, we obtain one well-known property of $n$-EP operators (see for example \cite[Theorem 3.3]{wangchun}).
\begin{corol}
  Let $T \in\mathcal{ R(H)}$   and let $n\ge 1$ be an integer. If $T$ is $n$-EP,
then $T$ is $(n + k)$-EP for an arbitrary integer $k\ge 0$.
\end{corol}

It well known that $T$ is EP operator if and only if $TT^{\dag 2}T=T^\dag T$ and $T^\dag T^{2}T^\dag=TT^\dag$, or equivalently $T$ and $T^{\dagger}$ are both hypo-EP operators. For class ($n$-EP) we have $T\in(n\text{-EP})$ if and only if $TT^{\dag }(T^\dag T^n)^*=(T^\dag T^n)^*$ and $T^\dag T^{n+1}T^\dag=T^{n}T^\dag$. For class of polynomially EP operators we have the following general setting.
%In next theorem provides further characterizations of polynomially EP operators.
\begin{theo}\label{equivEP}
Let $T\in \mathcal{R(H)}$ and $p\in Poly_{1}$ such that $p(0)=0.$ The  following assertions are equivalent:
\begin{itemize}
\item[(1)]\: $T$ is $p$-EP;
\item[(2)]\: $TT^{\dag}(T^\dag p(T))^*=(T^\dag p(T))^*$ and $T^{\dag} Tp(T)T^\dag=p(T)T^\dag$;
\item[(3)]\: $$p(T)TT^\dag+\left(p(T)TT^\dag\right)^*=p(T)+\overline{p}(T^*) \text{ and } T^{\dag} p(T)T+\left(T^{\dagger}p(T)T\right)^*=p(T)+\overline{p}(T^*).$$
\end{itemize}
\end{theo}
\begin{proof}
(1)$\Longrightarrow$(2) By assumption and from \cite[Theorem 2.5]{dija} we have,
$$\begin{aligned}
     R((T^\dag p(T))^*)\subset R(T)& \mbox{ and }R(p(T)T^\dag)\subset R(T^{*})\\
          \Longrightarrow\:\:&  P_{R(T)} (T^\dag p(T))^*=(T^\dag p(T))^* \mbox{ and } P_{R(T^{*})} p(T)T^\dag=p(T)T^\dag\\
          \Longrightarrow \:\:& TT^{\dag} (T^\dag p(T))^*=(T^\dag p(T))^* \mbox{ and } T^{\dag} T p(T)T^\dag=p(T)T^\dag.
    \end{aligned}$$
(2)$\Longrightarrow$(1) We have $$TT^{\dag} (T^\dag p(T))^*=(T^\dag p(T))^* \mbox{ and } T^{\dag} Tp(T)T^\dag=p(T)T^\dag.$$
Multiplying right by $T^*$ and by $T$ respectively, we get $$P_{R(T)}\overline{p}(T^*)=TT^{\dag} \overline{p}(T^*)=\overline{p}(T^*) \mbox{ and } P_{R(T^{*})}p(T)=T^{\dag} T p(T)=p(T).$$ Thus, $R(p(T))\subset R(T^{*})$ and  $R(\overline{p}(T^*))\subset R(T).$ Therefore by \cite[Theorem 2.5]{dija}, $T$ is $p$-EP.\\
(1)$\Longrightarrow$ (3) Is derived from  \cite[Theorem 2.5(iv)]{dija}.\\ %Is obvious.\\
(3)$\Longrightarrow$(1) Multiplying $p(T) T T^{\dagger}+\left(p(T) T T^{\dagger}\right)^*=p(T)+\overline{p}(T^*)$ by $T T^{\dagger}$ from the right side, we get
$$p(T) T T^{\dagger}+T T^{\dagger}\overline{p}(T^*) T T^{\dagger}=p(T) T T^{\dagger}+\overline{p}(T^*) T T^{\dagger},$$ which gives $T T^{\dagger}\overline{p}(T^*)=\overline{p}(T^*) T T^{\dagger}=\overline{p}(T^*)$. Thus, $p(T)= p(T)TT^{\dag}.$ In the same way, we show that  $T^{\dagger} p(T)T+\left(T^{\dagger}  p(T)T\right)^*=p(T)+\overline{p}(T^*)$ implies $p(T)=T^{\dag}Tp(T)$. Therefore by \cite[Theorem 2.5(iv)]{dija}, $T$ is $p$-EP.
\end{proof}
\begin{prop}\label{rest}
  Let $T\in \mathcal{R(H)}$ and $p\in Poly_{1}$ such that $p(0)=0$ and $R(p(T))$ is closed.
  \begin{enumerate}
    \item If $p(T)\in (EP)$ then $T\in (p\text{-EP})$.
    \item If $T\in \text{(SD)}$ then $T\in (p\text{-EP})$ if and only if $p(T)\in \text{(EP)}$.
  \end{enumerate}
\end{prop}
\begin{proof}
  (1)  If $p(T)\in \text{(EP)}$, then  $R(p(T))=R(p(T)^{*})\subseteq R(T^{*})$ and  $R(p(T)^{*})\subset R(T)$, using \cite[Theorem 2]{emz3} and \cite[Theorem 2.5]{dija} we get $T\in (p\text{-EP)}$.\\
  (2)  Suppose that $T\in (p\text{-EP})\cap\text{(SD)}$ and let $A=p(T),$ $B=T^{\dagger}.$ We have $A^{*}AN(B^{*})=p(T)^{*}p(T)N(T^{\dagger*})=p(T)^{*}p(T)N(T)={0}\subseteq N(B^{*}).$ On the
other hand, for $x\in N(A)=N(p(T))$ and since $T$ is SD, we have $ABB^{*}x=p(T)T^{\dagger}T^{\dagger*}x =T^{\dagger}p(T)T^{\dagger*}x=T^{\dagger}T^{\dagger*}p(T)x=0.$ Thus $BB^{*}N(A)\subseteq N(A).$ So, by
\cite[Remark 1]{emz3}, we derive that $(p(T)T^{\dagger})^{\dagger}=T^{\dagger\dagger}p(T)^{\dagger}= Tp(T)^{\dagger}.$ Similarly, for $A=T^{\dagger}$
and $B=p(T),$ we have $A^{*}AN(B^{*})\textcolor{blue}{\subseteq} N(B^{*})$ and $BB^{*}N(A)\subseteq N(A).$ Thus,
$(T^{\dagger}p(T))^{\dagger}=(p(T))^{\dagger}T.$ Since by assumption $T$ is $p$-EP, we get $Tp(T)^{\dagger}=p(T)^{\dagger}T$
and hence $p(T)p(T)^{\dagger}=p(T)^{\dagger}p(T)$ or equivalently $p(T)$ is EP.

\end{proof}
In the special case where $p(t)=t^{n}$ in Propsition \ref{rest}, we obtain one well-known
property of $n$-EP operators (see \cite[Theorem 3.6]{emz3}).

\section{Polynomially hypo-EP operators}

At the start of this section, we introduce and study the concept of polynomially hypo-EP operators, which extends the notion of $n$-hypo-EP operators.
  \begin{definition}
    Let $T\in \mathcal{R(H)}$ and $p\in Poly_{1}$ be nontrivial.
If $T^{\dagger}Tp(T)=p(T)T^{\dagger}T$, then we say that $T$ is
$p$-hypo-EP (or $p$-HEP, for short). An operator
$S\in \mathcal{L(H)}$ is polynomially hypo-EP, if S is $q$-hypo-EP for some $q\in Poly_{1}.$
  \end{definition}
  It is clear that both $T$ and $T^{*}$ are in hypo-EP if and only if $T$ is EP. Also, for $n$-EP operators, $T$ and $T^{*}$ are simultaneously in ($n$-HEP) if and only if $T$ is in ($n$-EP). For $p$-EP operators, by \cite[Theorem 2.5]{dija} and Theorem \ref{pHEP} below, we have for $p\in Poly_{1}$  such that $p(0)=0,$ $T$ is in ($p$-EP) if and only $T$ is in ($p$-HEP) and $T^{*}$ is in $(\bar{p}$-HEP).\\
 From  \cite[Theorem 2.4]{itoh}, $T\in(HEP)$ if and only if $T=T^{\dagger} T^2$, or equivalently, $R(T)\subset R(T^*)$ (\cite{ding1}).   %We proposed a wider class of operators which includes hypo-EP and $n$-hypo-EP operators.
 Obviously, when $T$ is $p$-HEP for $p(t)=t$ ( or $p(t) = t^{n}$), then $T$ is hypo-EP ($n$-hypo-EP).\\
The following proposition, we prove that if $T$ is hypo-EP, then $T$ is $p$-HEP for an arbitrary $p\in Poly_{1}.$
\begin{prop}
  Let $T\in \mathcal{R(H)}$ and let $p\in Poly_{1}$. If $T\in(HEP)$, then $T\in(p\text{-HEP)}.$
\end{prop}
\begin{proof}
Suppose that $T\in\text{(HEP)}$ and $p(t)=\sum_{i=0}^{n}a_{i}t^{j}$. We have
  $$
  \begin{aligned}
     p(T)T^{\dagger}T&=(a_{n}T^{n}+\cdots+a_{1}T+a_{0}I)T^{\dagger}T \\
      &=a_{n}T^{n}T^{\dagger}T+\cdots+a_{1}TT^{\dagger}T+a_{0}T^{\dagger}T\\
      &=a_{n}T^{n}+\cdots+a_{1}T+a_{0}T^{\dagger}T\\
      &=a_{n}T^{\dagger}TT^{n}+\cdots+a_{1}T^{\dagger}TT+a_{0}T^{\dagger}T\\
      &=T^{\dagger}Tp(T).
  \end{aligned}
  $$
  Therefore, $T$ is $p$-HEP.

\end{proof}

The converse of this implication is generally false. In the following examples, we observe that a polynomially hypo-EP operator is not hypo-EP in general.
\begin{example}
  Let $T$ be the following matrix defined by
    $T=\begin{pmatrix}
           1 & 0 & 2 & 3 \\
           0 & 2 & 0 & 0 \\
           0 & 0 & 0 & 3 \\
           0 & 0 & 0 & 0 \\
         \end{pmatrix}
  $
on $\mathcal{H}=\mathbb{C}^{4}.$ Since
  $$T^{\dagger}=\begin{pmatrix}
     \frac{1}{5} & 0 & -\frac{1}{5}  & 0 \\
     0 & \frac{1}{2}  & 0 & 0 \\
     \frac{2}{5}  & 0 & -\frac{2}{5}  & 0 \\
     0 & 0 & \frac{1}{3}  & 0 \\
   \end{pmatrix},\:\:
  T^{2}=\begin{pmatrix}
    1 & 0 & 2 & 9 \\
    0 & 4 & 0 & 0 \\
    0 & 0 & 0 & 0 \\
    0 & 0 & 0 & 0 \\
  \end{pmatrix}
\text{ and }
T^{3}=\begin{pmatrix}
  1 & 0 & 2 & 9 \\
  0 & 8 & 0 & 0 \\
  0 & 0 & 0 & 0 \\
  0 & 0 & 0 & 0 \\
\end{pmatrix},\:$$
  we get,
  $T^{3}-T^{2}=
  \begin{pmatrix}
  0 & 0 & 0 & 0 \\
  0 & 4 & 0 & 0 \\
  0 & 0 & 0 & 0 \\
  0 & 0 & 0 & 0 \\
\end{pmatrix}$
and $$(T^{3}-T^{2})T^{\dagger}T=
\begin{pmatrix}
   0 & 0 & 0 & 0 \\
   0 & 4 & 0 & 0 \\
   0 & 0 & 0 & 0 \\
   0 & 0 & 0 & 0 \\
 \end{pmatrix}
  =T^{\dagger}T(T^{3}-T^{2}).$$
 Thus, $T$ is $p$-HEP where $p(t)=t^{3}-t^{2}.$ Since $$TT^{\dagger}-T^{\dagger}T=\begin{pmatrix}
                               \frac{4}{5} & 0 & -\frac{2}{5} & 0 \\
                               0 & 0 & 0 & 0 \\
                               -\frac{2}{5} & 0 & \frac{1}{5} & 0 \\
                               0 & 0 & 0 & -1 \\
                             \end{pmatrix}
 $$ is not negative operator, we conclude that $T$ is not hypo-EP. Therefore $T$ is $p$-HEP but is not hypo-EP.
\end{example}
\begin{example}Let $p(t)=t^{3}+4it^{2}-(1+4i)t+2,$ $t\in \mathbb{C},$ and $$T=\begin{pmatrix}
           1 & 1 \\
           0 & 0 \\
         \end{pmatrix}$$
  on $\mathcal{H}=\mathbb{C}^{2}.$
  Since $T^{2}=T$, we obtain $p(T)=T^{3}+4iT^{2}-(1+4i)T+2I=2I.$ Thus, $T$ is $p$-HEP. Using $$T^{\dagger}=\begin{pmatrix}
                                             \frac{1}{2} & 0 \\
                                             \frac{1}{2} & 0 \\
                                           \end{pmatrix}$$
  notice that $$T^{\dagger}T-TT^{\dagger}=\begin{pmatrix}
                                            \frac{1}{2} & \frac{1}{2} \\
                                            0 & 0 \\
                                          \end{pmatrix}
  -\begin{pmatrix}
     1 & 0 \\
     0 & 0 \\
   \end{pmatrix}
  =\begin{pmatrix}
     -\frac{1}{2} & \frac{1}{2} \\
     0 & 0 \\
   \end{pmatrix}
$$
is not positive operator. Hence $T$ is not hypo-EP.
\end{example}
The next theorem provides characterizations of polynomially hypo-EP operators.
\begin{theo}\label{pHEP}
Let $T\in \mathcal{R(H)},$ $p\in Poly_{1}$ such that p(0)=0. The following statements are equivalent:
    \begin{enumerate}
       \item[(i)] $T\in (p\text{-HEP})$;
       \item [(ii)] $p(T)=T^{\dagger}Tp(T)$;
       \item[(iii)] $R(p(T)) \subset R(T^*)$;
       \item[(iv)]  $T^{\dag}T p(T)T^\dag=p(T)T^\dag$;
       \item[(v)] $T^{\dagger}Tp(T)+(T^{\dagger}Tp(T))^* =p(T) +p(T)^{*}.$
    \end{enumerate}
\end{theo}
\begin{proof}
    (i)$\Longleftrightarrow$(ii) Suppose that $T$ is $p$-HEP and $p(t)=\sum_{i=1}^{n}a_{i}t^{j}$. Then, $$p(T)T^{\dagger}T=(a_{n}T^{n}+\cdots+a_{1}T) T^{\dagger}T =(a_{n}T^{n-1}+\cdots+a_{1})TT^{\dagger}T= p(T).$$
     (ii)$\Longleftrightarrow$(iii) Evident because $p(T)=P_{R(T^*)}p(T)\iff R(p(T))\subset R(T^*)$).\\

  (iii)$\Longrightarrow$(iv) By the previous equivalence, we get $T^\dag Tp(T)=p(T)$, thus,   $T^\dag Tp(T)T^\dag=p(T)T^\dag$.

(iv)$\Longrightarrow$(iii) We have $$T^{\dag} Tp(T)T^\dag=p(T)T^\dag.$$
Multiplying right by $T$, we get $$P_{R(T^*)}p(T)=T^{\dag}Tp(T)=p(T).$$ Therefore, $R(p(T))\subset R(T^{*})$ .  \\

(ii)$\Longrightarrow$ (v) is obvious.\\
(v)$\Longrightarrow$(ii) Multiplying $T^{\dagger}Tp(T) +\left(T^{\dagger}Tp(T)\right)^*=p(T)+\left(p(T)\right)^*$ by $T^{\dagger}T$ from the right side, we obtain that $$T^{\dagger}Tp(T)+p(T)^* T^{\dagger}T T^{\dagger}T =p(T)T^{\dagger}T+(p(T))^*T^{\dagger}T,$$ which gives $T^{\dagger}Tp(T)=p(T)$.
\end{proof}
In the particular case that $p(t)=t^n$ in Theorem \ref{nHEP}, we  obtain new characterizations of $n$-hypo-EP operators with a closed range.
\begin{corol}
   Let $T\in \mathcal{R(H)}$. The following statements are equivalent:
    \begin{enumerate}
       \item [(i)] $T\in (n\text{-HEP})$;
       \item [(ii)]$R(T^n) \subset R(T^*)$;
       \item [(iii)] $T^{\dag}T^{n+1}T^\dag=T^nT^\dag $;
       \item [(iv)]$T^{\dagger}T^{n+1}+(T^{\dagger}T^{n+1})^* =T^n +(T^n)^{*}.$
    \end{enumerate}
\end{corol}
\begin{prop}\label{eqHEP}
  Let $T\in \mathcal{R(H)},$ $p\in Poly_{1}$ such that p(0)=0. The following statements are equivalent:
  \begin{enumerate}
    \item [(i)] $T\in (p\text{-HEP)}$;
    \item[(ii)]$p(T)^{*}=p(T)^{*}T^{\dagger}T$;
    \item[(iii)]$\|p(T)^{*}x\|\leq k\|Tx\|$ for some $k\geq0$ and all $x\in \mathcal{H}$;
    \item[(iv)]$\|p(T)^{*}x\|\leq k\|\omega(T)x\|$ for some $k\geq0$ and all $x\in \mathcal{H}.$
  \end{enumerate}
\end{prop}

\begin{proof}
  The equivalences (i) $\Longleftrightarrow$ (iii) and (i) $\Longleftrightarrow$ (iv) follow from the fact that $R(T^\dag)=R(T^*)$ and the results of Douglas \cite{douglas}.
\end{proof}

The class of hypo-EP operators is closed under unitary equivalence and restrictions to  reducing subspaces, see \cite{vinoth}. For $p$-hypo-EP class, we have the following  properties.

\begin{prop}
  Let  $T\in \mathcal{R(H)}$ be  unitarily equivalent to some $S\in \mathcal{R(H)}$ and let $p\in Poly_{1}.$
  \begin{enumerate}
    \item If $M$ is a reducing subspace for $T$, then, $T_{|M}$ is $p$-HEP;
    \item $S$ is $p$-HEP if and only if $T$ is $p$-HEP.
  \end{enumerate}
\end{prop}
\begin{proof}
\begin{enumerate}
  \item Evident.
  \item Suppose that $S=UTU^{*}$ for some unitary operator $U$. Since $S^{\dagger}=UT^{\dagger}U$ and $p(S)=Up(T)U^{*},$ we get
   $$\begin{aligned}
     p(S) &= S^{\dagger}Sp(S)\iff Up(T)U^{*}=UT^{\dagger}U^{*} UTU^{*}Up(T)U^{*}\\
      & \iff Up(T)U^{*}=UT^{\dagger}Tp(T)U^{*}\\
      &\iff p(T)=T^{\dagger}Tp(T).
   \end{aligned}.
   $$
Thus, $S$ is $p$-HEP if and only if $T$ is $p$-HEP.
\end{enumerate}

\end{proof}
In the following theorem\textcolor{red}{,} we show that a $p$-HEP operator is $p_{k}$-HEP, where $p_{k}(t)=t^{k}(p(t)-p(0)),$ for arbitrary $k\in \mathbb{N}.$
\begin{theo}
  Let $T\in \mathcal{R(H)}$ and $p\in Poly_{1}.$ If $T$ is $p$-HEP,
then $T$ is $p_{k}$-HEP, where $p_{k}(t)=t^{k}(p(t)-p(0)),$ for every $k\in \mathbb{N}$.
\end{theo}
\begin{proof}
Notice that for every $k\in \mathbb{N}$, $p_k(0)=0$. Suppose that $T$ is $p$-HEP and $p(t)=\sum_{i=0}^{n}a_{i}t^{i},$ where $a_{0},a_{1},\ldots,a_{n}\in\mathbb{C}.$\\
  For $k=1,$ $p_{1}(t)=t(p(t)-p(0))=a_{n}t^{n+1}+\cdots+a_{1}t^{2}.$ By assumption
   $$a_{n}T^{n}+\cdots+a_{1}T+a_{0}I=T^{\dagger}T( a_{n}T^{n}+\cdots+a_{1}T+a_{0}I).$$
   Multiplying by $T$ from the right hand side, we get $$a_{n}T^{n+1}+\cdots+a_{1}T^{2}+a_{0}T=T^{\dagger}T(a_{n}T^{n+1}+\cdots+a_{1}T^{2}+a_{0}T).$$
Therefore,
$$p_{1}(T)=T^{\dagger}Tp_{1}(T).$$
So, by Theorem \ref{pHEP}, $T$ is $p_{1}$-HEP.\\
Assume that $T$ is $p_{k}$-HEP, where $$p_{k}(t)=t^{k}(p_{k}(t)-p(0))=a_{n}t^{n+k}+\cdots+a_{1}t^{k+1}.$$ Then $$a_{n}T^{n+k}+\cdots+a_{1}T^{k+1}=T^{\dagger}T(a_{n}T^{n+k}+\cdots+a_{1}T^{k+1}),$$
multiplying by $T$ from the right hand side, we get
$$p_{k+1}(T)=T^{\dagger}T(a_{n}T^{n+k+1}+\cdots+a_{1}T^{k+2})=T^{\dagger}Tp_{k+1}.$$
Hence, $T$ is $p_{k+1}$-HEP. Therefore, $T$ is $p_{k}$-HEP for all $k\in \mathbb{N}.$
\end{proof}
\begin{corol}
     $(n\text{-HEP}) \subseteq (n+k\text{-HEP})$ for all $k\in \mathbb{N}.$
\end{corol}
The next proposition is a generalization of \cite[Proposition 13]{emz3} for $n$-hypo-EP operator.
\begin{prop}\label{prop7}
  Let $T$ be in $\mathcal{R(H)}$ and $p\in Poly_{1}$  such that $p(0)=0.$ If one of the following conditions is satisfied:
  \begin{itemize}
    \item[(1)] $[T^{\dagger}T,p(T)+T^{\dagger}]=0,$
    \item[(2)] $[T^{\dagger}T,p(T)+T^{*}]=0,$
  \end{itemize}
  then, $T\in (p\text{-HEP}).$
\end{prop}
\begin{proof}
Let us first note that since $p(0)=0$, then, $$p(T)T^\dag T=TT^\dag p(T)=p(T) \text{ and } p(T)^*TT^\dag=T^\dag T p(T)^*=p(T)^*.$$  Suppose that (1) is verified and $p(t)=\sum_{i=0}^{n}a_{i}t^{j},$ where $a_{0},a_{1},\ldots,a_{n}\in\mathbb{C}.$ Then,
  $$\begin{aligned}
    0=&[T^{\dagger}T,p(T)+T^{\dagger}] \\
    =&T^{\dagger}Tp(T)+T^{\dagger}TT^{\dagger}-p(T)T^{\dagger}T-T^{\dagger2}T\\
    =&T^{\dagger}Tp(T)+T^{\dagger}-p(T)-T^{\dagger2}T.
  \end{aligned}
  $$
  Multiplying by $T^{\dagger}$ from the right hand side, we get $$T^{\dagger}Tp(T)T^{\dagger}=p(T)T^{\dagger}.$$
  Thus, by Theorem \ref{pHEP}, we conclude that $T\in(p\text{-HEP})$.\\
  Suppose that (2) is satisfied, then
  $$\begin{aligned}
    0=&[T^{\dagger}T,p(T)+T^{*}] \\
    =&T^{\dagger}Tp(T)+T^{\dagger}TT^{*}-p(T)T^{\dagger}T-T^{*}T^{\dagger}T\\
    %=& T^{\dagger}Tp(T)+T^{*}-\textcolor{red}{p(T)}-T^{*}T^{\dagger}T\\
   % =&T^{\dagger}Tp(T)+T^{*}-\textcolor{red}{p(T)}-T^{*}T^{\dagger}T\\
    =&T^{\dagger}Tp(T)+T^{*}-p(T)-T^{*}T^{\dagger}T.
  \end{aligned}
  $$
   Multiplying right by $T^{\dagger},$  $$T^{\dagger}Tp(T)T^{\dagger}=p(T)T^{\dagger}.$$
  Thus, by Theorem \ref{pHEP}, we derive that $T\in(p\text{-HEP})$.\\
\end{proof}
In the following proposition, we provide sufficient conditions for a polynomially hypo-EP operator to be a polynomially EP operator.
\begin{prop}\label{prop8}
Let $T\in\mathcal{R(H)}$ and $p\in Poly_{1}$ such that $p(0)=0$. Suppose that $T\in(p\text{-HEP})$. If $T$ satisfies one of the following statements,
\begin{itemize}
    \item[(1)] $[T, p(T)T^\dag]=0;$
    \item [(2)]  $[TT^\dag, p(T)+T^{\dagger}]=0;$
    \item [(3)]$[TT^\dag, p(T)+T^{*}]=0;$
\end{itemize}
then $T\in(p\text{-EP}).$
\end{prop}
\begin{proof}
(1) We have
\begin{eqnarray*}
[T ,p(T)T^\dag]=0&\iff& Tp(T)T^\dag-p(T)T^\dag T=0\\
  &\iff & p(T)TT^\dag=p(T)\\
   &\iff & TT^\dag (p(T))^{*}=(p(T))^{*}\\
   &\iff & R(p(T)^{*})\subset R(T)\\
   &\iff & T^*\in (p\text{-HEP}).
\end{eqnarray*}
Since by assumption $T\in(p\text{-HEP})$, we derive that $T\in(p\text{-EP})$.\\
(2)  We have
 \begin{eqnarray*}
  0&=& [TT^{\dagger},p(T)+T^{\dagger}]\\
  &=& TT^{\dagger} p(T)+TT^{\dagger} T^{\dagger}-p(T)TT^{\dagger} -T^{\dagger}TT^{\dagger}\\
  &=& p(T)+T T^{\dagger2}-p(T)TT^{\dagger}-T^{\dagger}
\end{eqnarray*}
Multiplying  by $TT^{\dagger}$ from the left hand side, we get, $$
 p(T)= p(T)TT^{\dagger}.$$ Since $T$ is assumed to be $p$-HEP, the result derives from \cite[Theorem 2.5]{dija}.\\
(3) Similarly, we get this assertion.

\end{proof}
The next theorem is generalization of \cite[Theorem 3.1]{vinoth} for hypo-EP operator. We use Riesz representation to give a necessary and sufficient condition for an operator to be polynomially HEP operator.
\begin{theorem}\label{nHEP}
  If $T\in\mathcal{ L(H)}$ and  $p\in Poly_{1}$ is such that $p(0)=0$, then $T$ is $p$-HEP if and only if for each $x\in \mathcal{H}$, there exists $ c>0$ such that
\begin{equation}\label{NHEP}
  |\left\langle p(T)x,y \right\rangle|\leq c\|Ty\|, \text{ for all }  y\in \mathcal{H}.
\end{equation}
\end{theorem}
\begin{proof}
  Assume $T$ is $p$-HEP. If $x\in N( p(T)),$ then the result is trivial. Let $x\in \mathcal{H}$ such that $ x\notin N( p(T)).$ Then $ p(T)x\in R( p(T))\subset R(T^{*}).$ Therefore there exists a non-zero $s\in \mathcal{H}$ such that $T^{*}s= p(T)x.$ Then for all $y\in \mathcal{H},$
  $$|\left\langle p(T)x,y \right\rangle|=|\left\langle T^{*}s,y \right\rangle|=|\left\langle s,Ty \right\rangle|\leq\|s\|\|Ty\|.$$ Taking $c=\|s\|,$ we get $|\left\langle p(T)x,y \right\rangle|\leq k\|Ty\|.$\\
  Conversely, assume that for each $x\in \mathcal{H},$ there exists $c>0$ such that $$|\left\langle p(T)x,y \right\rangle|\leq c\|Ty\|.$$
  for all $y\in \mathcal{H}$. Let $x\in \mathcal{H}$ be fixed.

  Then, for all $y\in \mathcal{H}$, we have $$c\|Ty\|\geq|\left\langle x,p(T)^{*}y \right\rangle|=|f_{x}(p(T)^{*}y)|$$ where we set $f_{x}(p(T)^{*}y)=\left\langle p(T)^{*}y,x \right\rangle.$
  Hence $$|(p(T)f^{*}_{x})^{*}y)|\leq c\|(T^{*})^{*}y\|,$$\:for some $c>0$, for all $y\in\mathcal{ H}$. By Douglas's theorem,
  $$p(T)f^{*}_{x}=T^{*}D, $$
  for some $D\in \mathcal{L}(\mathbb{C},\mathcal{H}).$ Taking adjoint on both sides gives $$f_{x}p(T)^{*}=g_{x}T,$$ where $g_{x}=D^{*} \in \mathcal{L}(\mathcal{H},\mathbb{C}).$ By Riesz representation theorem, there exists $x' \in \mathcal{H}$ such that $$g_{x}(Ts)=\left\langle Ts,x' \right\rangle$$ for all $s\in\mathcal{H}.$ Hence for $s\in\mathcal{H},$ $$f_{x}p(T)^{*}s=g_{x}Ts$$
  implies that $$\left\langle p(T)^{*}s,x \right\rangle=\left\langle Ts,x'\right\rangle.$$ Hence for each $x\in \mathcal{H}$, there exists $x'\in\mathcal{ H}$ such that $p(T)x=T^{*}x'.$ Thus, $R(p(T))\subset R(T^{*}),$ i.e, $T$ is $p$-HEP.
\end{proof}
Applying Theorem \ref{nHEP}, for $p(t)=t^{n}$ we get the following consequence which involve characterization of $n$-hypo-EP operators.
\begin{corol}
  Let $T\in\mathcal{ R(H)}.$ Then $T$ is $n$-hypo-EP if and only if for each $x\in \mathcal{H}$, there exists $c>0$ such that
  \begin{equation}
  |\left\langle T^{n}x,y \right\rangle|\leq c\|Ty\|, \text{ for all }  y\in \mathcal{H}.
\end{equation}
\end{corol}
Recall that the closed range operator $T \in \mathcal{L(H)}$ has the following matrix representation
\begin{equation}\label{1}
 T=\left[\begin{array}{cc}
T_1 & 0 \\
T_2 & 0
\end{array}\right]:\left[\begin{array}{c}
R(T^*)  \\
N(T)
\end{array}\right] \rightarrow \left[\begin{array}{c}
R(T^*)  \\
N(T)
\end{array}\right],
\end{equation}
with respect to the orthogonal decomposition $$\mathcal{H}=R(T^*) \oplus N(T)$$
and $D=T_1^*T_1 +T_2^*T_2  \in \mathcal{L}(R(T^*))$ is an invertible operator. %In this case, the Moore-Penrose inverse of $T$ is represented by \textcolor{green}{ Oui vous avez raison, J ai  changé}
An easy computation shows that the Moore-Penrose inverse of $T$ is represented by
\begin{equation}\label{2}
T^{\dagger}=\left[\begin{array}{ll}
D^{-1}T_1^* &  D^{-1}T_2^*\\
0 & 0
\end{array}\right] .
\end{equation}

Using the representation (\ref{1}), we obtain the following characterizations of $p$-HEP operators.
\begin{theorem}\label{rep}
Let $T \in \mathcal{L(H)}$ be with a closed range. Suppose that $T$ is  represented by (\ref{1}),  $p(z)=a_n z^n+$ $a_{n-1} z^{n-1}+\cdots+a_0$ and $q(z)=a_n z^{n-1}+a_{n-1} z^{n-2}+\cdots+a_1$, where $a_0, a_1, \ldots, a_n \in$ $\mathbb{C}$. Then the  following statements are equivalent:
\begin{itemize}
    \item[(i)] $T$ is $p$-HEP;
    \item [(ii)]  $T_2 q(T_{1})=0$.
    \end{itemize}
    In addition, if $T$ is $p$-HEP then
$$ p(T)=\left[\begin{array}{cc}
p(T_1) & 0 \\
0 & a_{0}I
\end{array}\right].
$$
\end{theorem}
\begin{proof}
Since $T$ is assumed to be represented as in (\ref{1}), using the equality (\ref{2}) and the following representation
 $$ p(T)=\left[\begin{array}{cc}
p(T_1) & 0 \\
T_2 q(T_{1}) & a_{0}I
\end{array}\right],
$$
we derive that $$
\begin{aligned}
    T^{\dag}Tp(T)= &\left[\begin{array}{cc}
I & 0 \\
0 & 0
\end{array}\right]\left[\begin{array}{cc}
p(T_1) & 0 \\
T_2 q(T_{1}) & a_{0}I
\end{array}\right] \\
    =& \left[\begin{array}{cc}
p(T_1) & 0 \\
0 & 0
\end{array}\right].
\end{aligned}
$$
Furthermore, we have
 $$
\begin{aligned}
    p(T)T^{\dag}T= &\left[\begin{array}{cc}
p(T_1) & 0 \\
T_2 q(T_{1}) & a_{0}
\end{array}\right]\left[\begin{array}{cc}
I & 0 \\
0 & 0
\end{array}\right] \\
    =& \left[\begin{array}{cc}
p(T_1) & 0 \\
T_2 q(T_{1}) & 0
\end{array}\right].
\end{aligned}
$$
Hence, $ T^{\dag}Tp(T)=p(T)T^{\dag}T$ if and only if $T_2 q(T_{1})=0.$
\end{proof}
\begin{theorem}\label{equi}
    Let $T\in \mathcal{R(H)}$ and polynomials $p$, $q$ be as in Theorem \ref{rep}. The following statements are equivalent:
\begin{itemize}
    \item[(i)] $T$ is $p$-HEP;
    \item[(ii)] $T^{\dag}Tp(T)T^*T=p(T)T^*T;$
    \item[(iii)] $T^{\dag}Tp(T)T^*=p(T)T^*.$
\end{itemize}
\end{theorem}
\begin{proof}
(i)$\iff $(ii):  Let the polynomial $q$ be as in Theorem \ref{rep}.  Since
$$
\begin{aligned}
p(T)T^*T= &\left[\begin{array}{cc}
p(T_1) & 0 \\
T_2 q(T_{1}) & a_{0}I
\end{array}\right]\left[\begin{array}{cc}
T^{*}_1 & T^{*}_2\\
0 & 0
\end{array}\right]  \left[\begin{array}{cc}
T_1 & 0 \\
T_2 & 0
\end{array}\right]  \\
    =& \left[\begin{array}{cc}
p(T_1)T^{*}_1 & p(T_1)T^{*}_2 \\
T_2 q(T_{1})T^{*}_1 & T_2 q(T_{1})T^{*}_2
\end{array}\right]\left[\begin{array}{cc}
T_1 & 0 \\
T_2 & 0
\end{array}\right]\\
=&\left[\begin{array}{cc}
p(T_1)T^{*}_1T_1+p(T_1)T^{*}_2T_2 & 0 \\
T_2 q(T_{1})T^{*}_1T_1 +T_2 q(T_{1})T^{*}_2 T_2 & 0
\end{array}\right]\\
=&\left[\begin{array}{cc}
p(T_1)D & 0 \\
T_2 q(T_{1})D & 0
\end{array}\right].
\end{aligned}
    $$
And
$$
   T^{\dag}Tp(T)T^*T =\left[\begin{array}{cc}
I & 0 \\
0 & 0
\end{array}\right]\left[\begin{array}{cc}
p(T_1)D & 0 \\
T_2 q(T_{1})D & 0
\end{array}\right]
    = \left[\begin{array}{cc}
p(T_1)D & 0 \\
0 & 0
\end{array}\right].
$$
We deduce that  $T^{\dag}Tp(T)T^*T=p(T)T^*T$ if and only if $T_2 q(T_{1}) = 0$. Therefore by Theorem \ref{rep} we get equivalent.\\
(ii)$\iff $(iii): Multiplying the identities of assertions (ii) and (iii) from the right hand sides by $T$ and $T^\dag$ respectively allows us to find this equivalence.\\
\end{proof}
Applying Theorem \ref{equi}, we get the following consequences which involve characterizations of $p$-EP operators

\begin{corol}
Let $T\in\mathcal{R(H)}$ and $p\in Poly_{1}$. The following statements are equivalent:
\begin{itemize}
 \item[(i)] $T$ is $p$-EP;
    \item[(ii)] $T^{\dag}Tp(T)T^*T=p(T)T^*T$ and  $TT^{\dag}\overline{p}(T^*)TT^*=\overline{p}(T^*)TT^* ;$
    \item[(iii)] $T^{\dag}Tp(T)T^*=p(T)T^*$ and $TT^{\dag}\overline{p}(T^*)T=\overline{p}(T^*)T.$
\end{itemize}
\end{corol}
Theorem \ref{equi} implies new characterizations of $n$-hypo-EP and hypo-EP operators.
\begin{corol}
Let $T\in \mathcal{R(H)}$ and $n\geq1.$ The following statements are equivalent:
    \begin{itemize}
    \item[(i)] $T$ is $n$-hypo-EP;
    \item[(ii)] $T^{\dag}T^{n+1}T^*T=T^{n}T^*T;$
    \item[(iii)] $T^{\dag}T^{n+1}T^*=T^{n}T^*.$
\end{itemize}
In particular, the following conditions are equivalent:
    \begin{itemize}
    \item[(i)] $T$ is hypo-EP;
    \item[(ii)] $T^{\dag}T^2T^*T=TT^*T;$
    \item[(iii)] $T^{\dag}T^2T^*=TT^*.$
\end{itemize}
\end{corol}


\begin{thebibliography}{9}
\bibitem{Bourdon} {Bourdon, Paul, et al}, {\it{Closed-range posinormal operators and their products}} Linear Algebra and its Applications, \textbf{671}  (2023), 38-58.
 \bibitem{camb} {S. L. Campbell and C. D. Meyer}, {\it{EP Operators and Generalized Inverses}}, Canadian Mathematical Bulletin,  \textbf{18} (1975), 327-333.
\bibitem{chen}{S. Cheng, Y.Tian}, {\it{Two sets of new characterizations for normal and EP matrices}}.
Linear Algebr. Appl.\textbf{375}(2003), 181–195.
\bibitem{cvet}{M.D. Cvetkovi$\acute{c}$,  D. Mosi$\acute{c}$}, {\it{Drazin invertibility, characterizations and structure of polynomially normal operators.}}, {Linear Multilinear Algebr. } (2021).
 \bibitem{djor}{D.S. Djordjevi$\acute{c}$, M. Ch$\bar{o}$,  D. Mosi$\acute{c}$}, {\it{Polynomially normal operators}}, {Ann. Funct. Anal. } \textbf{11} (2020), 493–504.
\bibitem{douglas} {R.G. Douglas}, {\it{On majorization, factorization and range inclusion of operators on Hilbert space}}, {Proc. Amer. Math. Soc.}  \textbf{17.2} (1966), 413-415.
\bibitem{driv}{D. Drivaliaris, S. Karanasios, D. Pappas}, {\it{Factorizations of EP operators}}, {Linear Algebr. Appl.} \textbf{429(7)} (2008), {1555–1567}.
 \bibitem {emz1} { H. Ezzahraoui, M. Mbekhta and E. H. Zerouali}, {\it{On the Cauchy dual of closed range operators}}, {Acta Sci.Math},  \textbf{85} (2019), {231–248}.
 \bibitem {ferr} {D.E. Ferreyra, F.E. Levis, Saroj B. Malik, A.N. Priori}, {\it{On star-dagger matrices and the core-EP decomposition}}, Journal of Computational and Applied Mathematics, \textbf{ 438} (2024), 0377-0427.
  \bibitem{hartwig} {R. Hartwig and K. Spindelböck}, {\it{Matrices for which $A^*$ and $A^\dag$ commute}}, Linear and Multilinear Algebra, \textbf{14} (1983), 241-256.
\bibitem{menkad} {S. Menkad, and E. Anissa }. {\it{Some Results of n-EP Operators on Hilbert Spaces.}} International Journal of Analysis and Applications, \textbf{22} (2024), 78-78.
\bibitem{dija} {D.Mosić, and Miloš D. Cvetković}, {\it{Polynomially EP operators}}, Aequationes mathematicae \textbf{96.5} (2022): 1075-1087.
  \bibitem{itoh} {M. Itoh}. {\it{On some EP operators}}, {Nihonkai Mathematical Journal}, \textbf{16} (2005), 49-56.
\bibitem{jibril} {A. A. S. Jibril}, {\it{On $n$-Power Normal Operators}}, The Journal for Science and Engenering, \textbf{33} (2008), 247--253.
\bibitem{joh} {Johnson, P. Sam and Vinoth, A.}, {\it{Product and factorization of hypo-EP operators}} Special Matrices, vol. 6, no. 1, 2018, pp. 376-382.
\bibitem{katzand1} {M. H. Pearl}, {\it{On generalized inverses of matrices}}, Mathematical Proceedings of the Cambridge Philosophical Society, \textbf{62} (1966), 673-677.
  \bibitem{ding1} {A. B. Patel and M. P. Shekhawat}, {\it{Hypo-EP Operators}}, Indian Journal of Pure and Applied Mathematics,  \textbf{47} (2016), 73-84.
\bibitem{SLN}  {Saroj {B. Malik}, Laura Rueda and Néstor Thome}  {\it{The class of $m$-EP and $m$-normal matrices}}. Linear and Multilinear Algebra, \textbf{64.11} (2016) 2119-2132.
  \bibitem{schwerdtfeger} {H. Schwerdtfeger}, {\it{Introduction to Linear Algebra and the Theory of Matrices}}, P. Noordhoff N.V., Groningen-Holland, 1961.
  \bibitem {emz3} {R. Semmami, H. Ezzahraoui and E. H. Zerouali},   {\it{on some properties of $n$-EP and $n$-hypo-EP operators,}} arXiv, arXiv:2410.20920.
\bibitem{wangchun} {Wang Xiunan and Chunyuan Deng},  {\it{Properties of m-EP operators}}. Linear and Multilinear Algebra, \textbf{65.7} (2017) 1349-1361.
  \bibitem{vinoth}  {A. Vinoth and P. Sam Johnson}, {\it{On Sum and Restriction of Hypo-EP Operators}}, {Functional Analysis, Approximation and Computation} \textbf{9.1} (2017), 37–41.
\end{thebibliography}
\end{document}